\documentclass{amsart}
\usepackage{amsfonts}
\usepackage[pdftex]{graphicx}

\newtheorem{theorem}{Theorem}
\theoremstyle{plain}

\newtheorem{corollary}{Corollary}

\newtheorem{example}{Example}

\numberwithin{equation}{section}

\input{tcilatex}

\begin{document}
\title[Spherical Indicatrices of Involute of a Space Curve in Euclidean
3-Space]{Spherical Indicatrices of Involute of a Space Curve in Euclidean
3-Space}
\author{Y\i lmaz Tun\c{c}er}
\address[Y\i lmaz Tun\c{c}er, Serpil \"{U}nal and Murat Kemal Karacan]{Usak
University, Science and Arts Faculty\\
Department of Mathematics}
\email[Corresponding Author]{yilmaz.tuncer@usak.edu.tr}
\author{Serpil \"{U}nal}
\email[Serpil \"{U}nal]{atwo@atwoinst.edu}
\author{Murat Kemal Karacan}
\date{August 26, 2011}
\thanks{This paper is in final form and no version of it will be submitted
for publication elsewhere.}

\begin{abstract}
In this work, we studied the properties of the spherical indicatrices of
involute curve of a space curve and presented some characteristic properties
in the cases that involute curve and evolute curve are slant helices and
helices, spherical indicatrices are slant helices and helices and we
introduced new representations of spherical indicatrices.
\end{abstract}

\maketitle

\section{Introduction}

The specific curve pairs are the most popular subjects in curve and surface
theory and involute-evolute pair is one of them. We can see in most
textbooks various applications not only in curve theory but also in surface
theory and mechanic.

In this study, the spherical indicatrices of involute of a space curve are
given. In order to make a involute of a space curve and its evolute curve
slant helix, the feature that spherical indicatrices curve's need to have
are examined.

Let $\gamma :I\longrightarrow IR^{3}$ be a curve with $\gamma ^{\prime
}\left( s\right) \neq 0$, where $\gamma ^{\prime }\left( s\right) =d$ $%
\gamma \left( s\right) /ds$. The arc-lenght $s$ of a curve $\gamma \left(
s\right) $ is determined such that $\left\Vert \gamma ^{\prime }\left(
s\right) \right\Vert =1.$ Let us denote $T\left( s\right) =\gamma ^{\prime
}\left( s\right) $ and we call $T\left( s\right) $ a tangent vector of $%
\gamma $ at $\gamma (s)$. We define the curvature of $\gamma $ by $\kappa
\left( s\right) =\left\Vert \gamma ^{\prime \prime }\left( s\right)
\right\Vert $. If $\kappa \left( s\right) \neq 0,$ then the unit principal
normal vector $N\left( s\right) $ of the curve at $\gamma (s)$ is given by $%
\gamma ^{\prime \prime }\left( s\right) =\kappa \left( s\right) N\left(
s\right) $. The unit vector $B\left( s\right) =T\left( s\right) \Lambda
N\left( s\right) $ is called the unit binormal vector of $\gamma $ at $%
\gamma (s)$. Then we have the Frenet-Serret formulae 
\begin{equation*}
T^{\prime }=\kappa N\text{ , \ }N^{\prime }=-\kappa T+\tau B\text{ , \ }%
B^{\prime }=-\tau N
\end{equation*}%
where $\tau \left( s\right) $ is the torsion of $\gamma $ at $\gamma (s)$ 
\cite{5}.

The curve $\gamma $ is called evolute of $\widetilde{\gamma }$ if the
tangent vectors are orthogonal at the corresponding points for each $s\in
I\subset IR.$ In this case, $\widetilde{\gamma }$ is called involute of the
curve $\gamma $ and there exists a relationship between the position vectors
as%
\begin{equation*}
\widetilde{\gamma }(s^{\ast })=\gamma (s)+\lambda T(s)
\end{equation*}%
where $\lambda $ is the distance between the curves $\gamma $ and $%
\widetilde{\gamma }$ at the corresponding points for each $s.$ The pair of ($%
\widetilde{\gamma }$, $\gamma $) is called a involute-evolute pair. $\lambda 
$ is not a constant for involute-evolute pairs\cite{5}.

On the other hand, Izumiya and Takeuchi have introduced the concept of slant
helix by saying that the normal lines make a constant angle with a fixed
straight line. They characterize a slant helix if and only if the geodesic
curvature of the principal image of the principal normal indicatrix%
\begin{equation}
\Gamma =\frac{f^{\prime }}{\kappa \left( 1+f^{2}\right) ^{3/2}}  \label{1}
\end{equation}%
is a constant function, where $f=\frac{\tau }{\kappa }$\cite{4,6} .

In this study, we denote $T$, $N$, $B,$ $\kappa $, $\tau $ and $\widetilde{T}
$, $\widetilde{N}$, $\widetilde{B},$ $\widetilde{\kappa }$, $\widetilde{\tau 
}$ are the Frenet equipments of $\gamma $ and $\widetilde{\gamma },$
respectively. Tangent, principal normal and binormal vectors are described
for the spherical curves which are called tangent, principal normal and
binormal indicatrices both the curves $\gamma $ and $\widetilde{\gamma },$
respectively. Throughout this study, both involute and evolute curves are
regular.

\section{Spherical indicatrices of involute of a curve}

In this section, we introduced the spherical indicatrices of involute curve
of a curve in Euclidean 3-space and gave considerable results by using the
properties of the curves, similar to the previous section. Let $\gamma $ be
a curve with its involute curve $\widetilde{\gamma }$ then 
\begin{equation}
\widetilde{\gamma }(s^{\ast })=\gamma (s)+\lambda T(s)  \label{2}
\end{equation}%
where 
\begin{equation}
\lambda =\left\vert c-s\right\vert  \label{3}
\end{equation}%
and $\left( c-s\right) $ is definitely positive. Let $\epsilon $ be the sign
of $\left( c-s\right) $ such that if $c-s>0$, $\epsilon =+1$ and if $c-s<0$, 
$\epsilon =-1.$ We differentiate the equation (\ref{2}) with respect to $s$,
we get%
\begin{equation*}
\widetilde{T}\frac{ds^{\ast }}{ds}=\left( 1-\epsilon \right) T+\kappa
\epsilon \left( c-s\right) N.
\end{equation*}%
Since $T$ and $\widetilde{T}$ are orthogonal, there is no any component of $%
\widetilde{T}$ on $T$. Thus $\epsilon $ has to be $+1$.

\begin{theorem}
\label{t1}Let $\widetilde{\gamma }$ be involute of a space curve, then we
have Frenet formula:%
\begin{equation*}
\widetilde{T}^{\prime }=\widetilde{\kappa }\widetilde{N},\text{ \ \ }%
\widetilde{N}^{\prime }=-\widetilde{\kappa }\widetilde{T}+\widetilde{\tau }%
\widetilde{B},\text{ \ \ }\widetilde{B}^{\prime }=-\widetilde{\tau }%
\widetilde{N}
\end{equation*}%
where 
\begin{equation*}
\widetilde{T}=N,\text{ \ }\ \widetilde{N}=\frac{-1}{\sqrt{1+f^{2}}}T+\frac{f%
}{\sqrt{1+f^{2}}}B,\ \ \ \widetilde{B}=\frac{f}{\sqrt{1+f^{2}}}T+\frac{1}{%
\sqrt{1+f^{2}}}B
\end{equation*}%
with the parametrization\qquad \qquad \qquad \qquad \qquad \qquad\ 
\begin{equation}
\frac{ds}{ds^{\ast }}=\frac{1}{\kappa \left( c-s\right) }\ \ \   \label{4}
\end{equation}%
$\ $and the curvature and torsion of $\widetilde{\gamma }\left( s^{\ast
}\right) $ are 
\begin{equation}
\widetilde{\kappa }=\frac{\sqrt{1+f^{2}}}{c-s},\text{ \ \ }\widetilde{\tau }=%
\frac{f^{\prime }}{\kappa \left( c-s\right) \left( 1+f^{2}\right) }.
\label{5}
\end{equation}%
The geodesic curvature of the the principal image of the principal normal
indicatrix of involute curve is%
\begin{equation}
\widetilde{\Gamma }=\frac{\left\{ \left( 1+f^{2}\right) \left( f^{\prime
\prime }\kappa -f^{\prime }\kappa ^{\prime }\right) -3\kappa f\left(
f^{^{\prime }}\right) ^{2}\right\} \left( 1+f^{2}\right) ^{\frac{3}{2}}}{%
\left( \kappa ^{2}\left( 1+f^{2}\right) ^{3}+\left( f^{^{\prime }}\right)
^{2}\right) ^{\frac{3}{2}}}  \label{6}
\end{equation}
\end{theorem}

From (\ref{5}), it is obvious that involute of $\gamma $ is a planar curve
if and only if $\gamma $ is a generalized helix. For further usage we denote 
$\frac{\text{\ }\widetilde{\tau }}{\widetilde{\kappa }}$ as $\widetilde{f}$.
By using (\ref{1}) and (\ref{5}) we obtained the relation%
\begin{equation*}
\widetilde{f}=\frac{f^{\prime }}{\kappa \left( 1+f^{2}\right) ^{3/2}}=\Gamma
\end{equation*}%
and so we have 
\begin{equation}
\widetilde{f}=\Gamma .  \label{6a}
\end{equation}%
Thus we have the following theorem.

\begin{theorem}
\label{c2}Let $\widetilde{\gamma }$ be involute of a space curve $\gamma $
then $\widetilde{\gamma }$ is a generalized helix if and only if its evolute
is a slant helix.
\end{theorem}

We obtained the relation between $\widetilde{\Gamma }$ and $\Gamma $\ by
using (\ref{1}) and (\ref{6}) as%
\begin{equation}
\widetilde{\Gamma }=\frac{\kappa ^{2}\left( 1+f^{2}\right) ^{4}\Gamma
^{\prime }}{\left( \kappa ^{2}\left( 1+f^{2}\right) ^{3}+\left( f^{\prime
}\right) ^{2}\right) ^{\frac{3}{2}}}  \label{7}
\end{equation}%
Thus we can give the following theorem.

\begin{theorem}
\label{c3}Let $\widetilde{\gamma }$ be involute of a space curve $\gamma $
then $\widetilde{\gamma }$ is a slant helix if and only if $\gamma $ is a
slant helix.
\end{theorem}

The spherical image of tangent indicatrix of $\widetilde{\gamma }$ is%
\begin{equation*}
\widetilde{\gamma }_{t}(s_{t}^{\ast })=\widetilde{T}(s^{\ast })
\end{equation*}%
with the natural parameter $s_{t}^{\ast }.$

\begin{theorem}
\label{t2}If the Frenet frame of the tangent indicatrix $\widetilde{\gamma }%
_{t}=\widetilde{T}$ \ of involute of $\gamma (s)$ is $\{\widetilde{T}_{t},%
\widetilde{N}_{t},\widetilde{B}_{t}\}$, we have Frenet formula:%
\begin{equation*}
\widetilde{T}_{t}^{\prime }=\widetilde{\kappa }_{t}\widetilde{N}_{t},\text{
\ \ }\widetilde{N}_{t}^{\prime }=-\widetilde{\kappa }_{t}\widetilde{T}_{t}+%
\widetilde{\tau }_{t}\widetilde{B}_{t},\text{ \ \ }\widetilde{B}_{t}^{\prime
}=-\widetilde{\tau }_{t}\widetilde{N}_{t}
\end{equation*}%
where%
\begin{equation}
\widetilde{T}_{t}=\widetilde{N}\text{ \ , \ }\widetilde{N}_{t}=\frac{-%
\widetilde{T}+\widetilde{f}\widetilde{B}}{\sqrt{1+\widetilde{f}^{2}}}\text{
\ , \ }\widetilde{B}_{t}=\frac{\widetilde{f}\widetilde{T}+\widetilde{B}}{%
\sqrt{1+\widetilde{f}^{2}}}  \label{8}
\end{equation}%
with the parametrization%
\begin{equation}
\frac{ds^{\ast }}{ds_{t}^{\ast }}=\frac{1}{\widetilde{\kappa }}  \label{9}
\end{equation}%
and the curvature and torsion of $\widetilde{\gamma }_{t}$ are%
\begin{equation}
\widetilde{\kappa }_{t}=\sqrt{1+\widetilde{f}^{2}},\text{\ }\widetilde{\tau }%
_{t}=\frac{\widetilde{f}^{\prime }}{\widetilde{\kappa }\left( 1+\widetilde{f}%
^{2}\right) }.  \label{10}
\end{equation}%
The geodesic curvature of the principal image of the principal normal
indicatrix of $\widetilde{\gamma }_{t}$ is%
\begin{equation}
\widetilde{\Gamma }_{t}=\frac{\left\{ \left( 1+\widetilde{f}^{2}\right)
\left( \widetilde{f}^{\prime \prime }\widetilde{\kappa }-\widetilde{f}%
^{\prime }\widetilde{\kappa }^{\prime }\right) -3\widetilde{\kappa }%
\widetilde{f}\left( \widetilde{f}^{^{\prime }}\right) ^{2}\right\} \left( 1+%
\widetilde{f}^{2}\right) ^{\frac{3}{2}}}{\left( \widetilde{\kappa }%
^{2}\left( 1+\widetilde{f}^{2}\right) ^{3}+\left( \widetilde{f}^{^{\prime
}}\right) ^{2}\right) ^{\frac{3}{2}}}  \label{11}
\end{equation}
\end{theorem}

Form (\ref{10}) we have the following theorem.

\begin{theorem}
\label{c5}Let $\widetilde{\gamma }$ be involute of a space curve $\gamma $
then spherical image of the tangent indicatrix of $\widetilde{\gamma }$ is a
circle on unit sphere if and only if evolute of $\widetilde{\gamma }$ is a
generalized helix.
\end{theorem}

By using (\ref{6a}) and theorem \ref{c5}, we can state the following
corollary.

\begin{corollary}
\label{c6}Spherical image of the tangent indicatrix of $\widetilde{\gamma }$
is a circle on unit sphere if and only if evolute of $\widetilde{\gamma }$
is a slant helix.
\end{corollary}

From (\ref{10}) and (\ref{11}) we obtained 
\begin{equation}
\frac{\widetilde{\tau }_{t}}{\widetilde{\kappa }_{t}}=\frac{\widetilde{f}%
^{^{\prime }}}{\widetilde{\kappa }\left( 1+\widetilde{f}^{2}\right) ^{1/2}}=%
\widetilde{\Gamma }  \label{11a}
\end{equation}%
and the relation between $\widetilde{\Gamma }$ and $\widetilde{\Gamma }_{t}$
is%
\begin{equation}
\widetilde{\Gamma }_{t}=\frac{\widetilde{\kappa }^{2}\left( 1+\widetilde{f}%
^{2}\right) ^{4}\widetilde{\Gamma }^{\prime }}{\left( \widetilde{\kappa }%
^{2}\left( 1+\widetilde{f}^{2}\right) ^{3}+\left( \widetilde{f}^{^{\prime
}}\right) ^{2}\right) ^{\frac{3}{2}}}  \label{12}
\end{equation}%
so we have the following theorem.

\begin{theorem}
Let $\widetilde{\gamma }$ be involute of a space curve $\gamma $ then
spherical image of the tangent indicatrix of $\widetilde{\gamma }$ is a
spherical helix if and only if involute of $\gamma $ is a slant helix. In
this case, spherical image of the tangent indicatrix of $\widetilde{\gamma }$
is a slant helix on unit sphere too.
\end{theorem}

The relation 
\begin{equation*}
\frac{\widetilde{\tau }_{t}}{\widetilde{\kappa }_{t}}=\frac{\kappa
^{2}\left( 1+f^{2}\right) ^{4}\Gamma ^{\prime }}{\left( \kappa ^{2}\left(
1+f^{2}\right) ^{3}+\left( f^{\prime }\right) ^{2}\right) ^{\frac{3}{2}}}
\end{equation*}%
can obtained by using (\ref{7}) and (\ref{11a}). Hence we can give the
following theorem.

\begin{theorem}
Let $\widetilde{\gamma }$ be involute of a space curve $\gamma $ then
followings are true.
\end{theorem}

\begin{theorem}
\textbf{i.} If $\gamma $ is a slant helix then spherical image of the
tangent indicatrix of its involute is a generalized helix.

\textbf{ii.} If $\widetilde{\gamma }$ is a slant helix then its spherical
image of the tangent indicatrix is a spherical slant helix.
\end{theorem}

The spherical image of principal normal indicatrix of involute of the curve $%
\gamma (s)$ is%
\begin{equation*}
\widetilde{\gamma }_{n}\left( s_{n}^{\ast }\right) =\widetilde{N}(s^{\ast })
\end{equation*}%
with the natural parameter $s_{n}^{\ast }.$

\begin{theorem}
If the Frenet frame of the principal normal indicatrix $\widetilde{\gamma }%
_{n}=\widetilde{N}$ of involute of the curve $\gamma (s)$ is $\{\widetilde{T}%
_{n},\widetilde{N}_{n},\widetilde{B}_{n}\}$, we have Frenet formula:%
\begin{equation}
\widetilde{T}_{n}^{\prime }=\widetilde{\kappa }_{n}\widetilde{N}_{n},\text{
\ \ }\widetilde{N}_{n}^{\prime }=-\widetilde{\kappa }_{n}\widetilde{T}_{n}+%
\widetilde{\tau }_{n}\widetilde{B}_{n},\text{ \ \ }\widetilde{B}_{n}^{\prime
}=-\widetilde{\tau }_{n}\widetilde{N}_{n}  \label{13}
\end{equation}%
where%
\begin{eqnarray}
\widetilde{T}_{n} &=&\frac{-\widetilde{T}+\widetilde{f}\widetilde{B}}{\sqrt{%
1+\widetilde{f}^{2}}}  \notag \\
\widetilde{N}_{n} &=&\frac{\widetilde{\kappa }\widetilde{f}\widetilde{f}%
^{^{\prime }}\left( 1+\widetilde{f}^{2}\right) }{\rho }\widetilde{T}-\frac{%
\widetilde{\kappa }^{2}\left( 1+\widetilde{f}^{2}\right) ^{3}}{\rho }%
\widetilde{N}+\frac{\widetilde{\kappa }\widetilde{f}^{^{\prime }}\left( 1+%
\widetilde{f}^{2}\right) }{\rho }\widetilde{B}  \label{14} \\
\widetilde{B}_{n} &=&\frac{\widetilde{\kappa }^{2}\widetilde{f}\left( 1+%
\widetilde{f}^{2}\right) ^{\frac{5}{2}}}{\rho }\widetilde{T}+\frac{%
\widetilde{\kappa }\widetilde{f}^{^{\prime }}\left( 1+\widetilde{f}%
^{2}\right) ^{\frac{3}{2}}}{\rho }\widetilde{N}+\frac{\widetilde{\kappa }%
^{2}\left( 1+\widetilde{f}^{2}\right) ^{\frac{5}{2}}}{\rho }\widetilde{B} 
\notag
\end{eqnarray}%
with the parametrization%
\begin{equation}
\frac{ds^{\ast }}{ds_{n}^{\ast }}=\frac{1}{\widetilde{\kappa }\sqrt{1+%
\widetilde{f}^{2}}}  \label{15}
\end{equation}%
and the curvature and torsion of $\widetilde{\gamma }_{n}$ are%
\begin{eqnarray}
\widetilde{\kappa }_{n} &=&\frac{\rho }{\widetilde{\kappa }^{2}\left( 1+%
\widetilde{f}^{2}\right) ^{3}}  \label{16} \\
\widetilde{\tau }_{n} &=&\left\{ \left[ \tfrac{-\widetilde{f}\left( 
\widetilde{f}^{^{\prime }}\right) ^{2}\widetilde{\kappa }^{2}\left( 1+%
\widetilde{f}^{2}\right) ^{3}}{\rho ^{2}}\right] +\left[ \left( \tfrac{-%
\widetilde{f}^{^{\prime }}\left( 1+\widetilde{f}^{2}\right) ^{\frac{3}{2}}}{%
\rho }\right) \left( \tfrac{\widetilde{\kappa }^{2}\left( 1+\widetilde{f}%
^{2}\right) ^{\frac{5}{2}}}{\rho }\right) ^{^{\prime }}\right] +\left[
\left( \tfrac{\widetilde{\kappa }\widetilde{f}^{^{\prime }}\left( 1+%
\widetilde{f}^{2}\right) ^{\frac{3}{2}}}{\rho }\right) ^{^{\prime }}\left( 
\tfrac{\widetilde{\kappa }\left( 1+\widetilde{f}^{2}\right) ^{\frac{5}{2}}}{%
\rho }\right) \right] \right\}  \notag
\end{eqnarray}
where 
\begin{equation*}
\rho =\sqrt{\left( \widetilde{f}^{^{\prime }}\right) ^{2}+\widetilde{\kappa }%
^{2}\left( 1+\widetilde{f}^{2}\right) ^{3}}
\end{equation*}%
$.$The geodesic curvature of the principal image of the principal normal
indicatrix of $\widetilde{\gamma }_{n}$ is%
\begin{equation}
\widetilde{\Gamma }_{n}=\frac{\widetilde{\kappa }^{2}\left\{ \left( 1+%
\widetilde{f}^{2}\right) \left( \sigma ^{\prime }\widetilde{\kappa }+2\sigma 
\widetilde{\kappa }^{\prime }\right) +6\sigma \widetilde{\kappa }\widetilde{f%
}\widetilde{f}^{^{\prime }}\right\} \left( 1+\widetilde{f}^{2}\right) ^{%
\frac{9}{2}}}{\left( \rho ^{2}+\sigma ^{2}\widetilde{\kappa }^{4}\left( 1+%
\widetilde{f}^{2}\right) ^{6}\right) ^{\frac{3}{2}}}  \label{17}
\end{equation}%
where
\end{theorem}

\begin{equation*}
\widetilde{\tau }_{n}=\sigma .
\end{equation*}

The spherical image of binormal indicatrix of involute of the curve $\gamma
(s)$ is%
\begin{equation*}
\widetilde{\gamma }_{b}\left( s_{b}^{\ast }\right) =\widetilde{B}\left(
s^{\ast }\right)
\end{equation*}%
with the natural parameter $s_{b}^{\ast }.$

\begin{theorem}
If the Frenet frame of the binormal indicatrix $\widetilde{\gamma }_{b}=%
\widetilde{B}$ of involute of the curve $\gamma (s)$ is $\{\widetilde{T}_{b},%
\widetilde{N}_{b},\widetilde{B}_{b}\}$, we have Frenet formula:%
\begin{equation*}
\widetilde{T}_{b}^{\prime }=\widetilde{\kappa }_{b}\widetilde{N}_{b},\text{
\ \ }\widetilde{N}_{b}^{\prime }=-\widetilde{\kappa }_{b}\widetilde{T}_{b}+%
\widetilde{\tau }_{b}\widetilde{B}_{b},\text{ \ \ }\widetilde{B}_{b}^{\prime
}=-\widetilde{\tau }_{b}\widetilde{N}_{b}
\end{equation*}%
where\ \ \ 
\begin{equation}
\widetilde{T}_{b}=-\widetilde{N}\text{ \ , \ }\widetilde{N}_{b}=\frac{-%
\widetilde{T}+\widetilde{f}\widetilde{B}}{\sqrt{1+\widetilde{f}^{2}}}\text{
\ , \ }\widetilde{B}_{b}=\frac{-\widetilde{f}\widetilde{T}-\widetilde{B}}{%
\sqrt{1+\widetilde{f}^{2}}}  \label{18}
\end{equation}%
with the parametrization\ \ \ \ \ \ \ \ \ 
\begin{equation}
\frac{ds^{\ast }}{ds_{b}^{\ast }}=\frac{1}{\widetilde{\tau }}  \label{19}
\end{equation}%
and the curvature and torsion of $\widetilde{\gamma }_{b}$ are%
\begin{equation}
\widetilde{\kappa }_{b}=\frac{\sqrt{1+\widetilde{f}^{2}}}{\widetilde{f}}\ ,\ 
\widetilde{\tau }_{b}=\frac{-\widetilde{f}^{^{\prime }}}{\widetilde{\kappa }%
\ \widetilde{f}\left( 1+\widetilde{f}^{2}\right) }.  \label{20b}
\end{equation}%
The geodesic curvature of the principal image of the principal normal
indicatrix of $\widetilde{\gamma }_{b}$ is%
\begin{equation}
\widetilde{\Gamma }_{b}=\frac{\left( 1+\widetilde{f}^{2}\right) ^{\frac{5}{2}%
}\left( -\widetilde{f}\left( \widetilde{f}^{^{\prime }}\right) ^{2}%
\widetilde{\tau }-\widetilde{\tau }^{2}\widetilde{f}^{2}\left( \frac{%
\widetilde{f}^{^{\prime }}}{\widetilde{\tau }}\right) ^{^{\prime }}\right) +3%
\widetilde{\tau }\widetilde{f}^{3}\left( \widetilde{f}^{^{\prime }}\right)
^{2}\left( 1+\widetilde{f}^{2}\right) ^{\frac{3}{2}}}{\left( \widetilde{\tau 
}^{2}\left( 1+\widetilde{f}^{2}\right) ^{3}+\widetilde{f}^{2}\left( 
\widetilde{f}^{^{\prime }}\right) ^{2}\right) ^{\frac{3}{2}}}.  \label{21}
\end{equation}
\end{theorem}

From (\ref{20b}) we have the following theorem.

\begin{theorem}
Let $\gamma $ be a space curve and $\widetilde{\gamma }$ be its involute
with nonzero torsion then spherical image of binormal indicatrix of $%
\widetilde{\gamma }$ is a circle on unit sphere if and only if $\widetilde{%
\gamma }$ is a generalized helix.
\end{theorem}

The relation 
\begin{equation}
\frac{\widetilde{\tau }_{b}}{\widetilde{\kappa }_{b}}=\frac{-\widetilde{f}%
^{^{\prime }}}{\widetilde{\kappa }\ \left( 1+\widetilde{f}^{2}\right) ^{3/2}}%
=-\widetilde{\Gamma }  \label{21a}
\end{equation}%
can easly obtaine from (\ref{20b}) and it gives us the following theorem.

\begin{theorem}
Let non-planar $\widetilde{\gamma }$ be involute of a space curve $\gamma $
then its spherical image of binormal indicatrix is spherical helix if and
only if evolute of $\gamma $ is a slant helix.
\end{theorem}

From (\ref{19}), (\ref{21}) and (\ref{21a}) we can obtaine the relation
between $\widetilde{\Gamma }$ and $\widetilde{\Gamma }_{b}$ which is 
\begin{equation}
\widetilde{\Gamma }_{b}=\frac{-\widetilde{\Gamma }^{\prime }}{\widetilde{%
\kappa }\sqrt{1+\widetilde{f}^{2}}\left( 1+\widetilde{\Gamma }^{2}\right)
^{3/2}}.  \label{22}
\end{equation}

Thus we have the following theorem.

\begin{theorem}
Spherical image of binormal indicatrix $\widetilde{\gamma }_{b}$ is a
spherical slant helix if and only if $\widetilde{\gamma }$ is a slant helix.
\end{theorem}

By using (\ref{4}), (\ref{5}), (\ref{6a}) and (\ref{20b}) we get the relation%
\begin{equation*}
\widetilde{\tau }_{b}=\frac{-\Gamma ^{^{\prime }}}{\ \kappa \Gamma \sqrt{%
1+f^{2}}\left( 1+\Gamma ^{2}\right) }
\end{equation*}%
and we have the following theorem.

\begin{theorem}
Spherical image of binormal indicatrix $\widetilde{\gamma }_{b}$ is circle
if and only if $\gamma $ is a slant helix.
\end{theorem}

\begin{corollary}
There are the following equations between the Frenet vectors of the
spherical indicatrices of an involute curve in Euclidean 3-space.
\end{corollary}

\begin{corollary}
\textbf{i.} $\widetilde{T}_{t}=-\widetilde{T}_{b},$

\textbf{ii.} $\widetilde{N}_{t}=\widetilde{T}_{n}=\widetilde{N}_{b}$

\textbf{iii.} \ $\widetilde{B}_{t}=-\widetilde{B}_{b}$
\end{corollary}

\begin{theorem}
\begin{theorem}
Let $\ \alpha $\ and $\ \beta \ $be two regular curves in $E^{3}.$Then $%
\alpha $\ and $\ \beta $ are similar curves with variable transformation if
and only if the principal normal vectors are the same for all curves\ \ \ \
\ \ \ \ \ \ \ \ \ \ \ \ \ \ \ \ \ \ \ \ \ \ \ \ \ \ \ \ \ \ \ \ \ \ \ \ \ \
\ \ \ \ \ \ \ \ \ \ \ \ \ \ \ \ \ \ \ \ \ \ \ \ \ \ \ \ \ \ \ \ \ \ \ \ \ \
\ \ \ \ \ \ 
\begin{equation*}
N_{\beta }\left( s_{\beta }\right) =N_{\alpha }\left( s_{\alpha }\right) \ \
\ \ \ \ \ \ \ \ 
\end{equation*}%
\ under the particular variable transformation\ \ \ \ \ \ \ \ \ \ \ \ \ \ \
\ \ \ \ \ \ \ \ \ \ \ \ \ \ \ \ \ \ \ \ \ \ \ \ \ \ \ \ \ \ \ \ \ \ \ \ \ \
\ \ \ \ \ \ \ \ \ \ \ \ \ \ \ \ \ \ \ \ \ \ \ \ \ \ \ \ \ $\ \ \ \ \ \ \ $%
\begin{equation*}
\frac{ds_{\beta }}{ds_{\alpha }}=\frac{\kappa _{\alpha }}{\kappa _{\beta }}\ 
\end{equation*}%
of the arc-lengths.
\end{theorem}
\end{theorem}

\ We can give the following corollary by using this theorem.

\begin{corollary}
Let $\ \widetilde{\gamma }_{t}$\ and $\ \widetilde{\gamma }_{b}\ $be two
regular curves in $E^{3}.$Then $\widetilde{\gamma }_{t}$\ and $\ \widetilde{%
\gamma }_{b}$ are similar curves with variable transformation if and only if
the principal normal vectors are the same for all curves\ \ \ \ \ \ \ \ \ \
\ \ \ \ \ \ \ \ \ \ \ \ \ \ \ \ \ \ \ \ \ \ \ \ \ \ \ \ \ \ \ \ \ \ \ \ \ \
\ \ \ \ \ \ \ \ \ \ \ \ \ \ \ \ \ \ \ \ \ \ \ \ \ \ \ \ \ \ \ \ \ \ \ \ \ \ 
\begin{equation*}
\ \widetilde{N}_{t}\left( s_{t}^{\ast }\right) =\widetilde{N}_{b}\left(
s_{b}^{\ast }\right) \ \ \ \ \ \ \ \ 
\end{equation*}%
\ under the particular variable transformation\ \ \ \ \ \ \ \ \ \ \ \ \ \ \
\ \ \ \ \ \ \ \ \ \ \ \ \ \ \ \ \ \ \ \ \ \ \ \ \ \ \ \ \ \ \ \ \ \ \ \ \ \
\ \ \ \ \ \ \ \ \ \ \ \ \ \ \ \ \ \ \ \ \ \ \ \ \ \ \ \ \ $\ \ \ \ \ \ \ $%
\begin{equation*}
\frac{ds_{b}^{\ast }}{ds_{t}^{\ast }}=\frac{\widetilde{\kappa }_{t}}{%
\widetilde{\kappa }_{b}}\ 
\end{equation*}%
of the arc-lengths.
\end{corollary}

\begin{example}
In \cite{4a}, the general equation of spherical helix family is obtained by
Monterde which is, 
\begin{equation*}
\alpha _{c}\left( s\right) =\left( \cos \left( s\right) \cos \left(
ws\right) +\frac{1}{w}\sin \left( s\right) \sin \left( ws\right) ,-\cos
\left( s\right) \sin \left( ws\right) +\frac{1}{w}\sin \left( s\right) \cos
\left( ws\right) ,\frac{1}{cw}\sin \left( s\right) \right)
\end{equation*}%
where $w=\frac{\sqrt{1+c^{2}}}{c},$ \ $c\in IR-\left\{ 0\right\} .$In \cite%
{4b},\ Kula et al. obtained the general equation of a slant helix family
similar to following%
\begin{equation*}
\gamma _{\mu }\left( s\right) =\left( \gamma _{\mu }^{1}(s),\gamma _{\mu
}^{2}(s),\gamma _{\mu }^{3}(s)\right)
\end{equation*}%
where 
\begin{eqnarray*}
\gamma _{\mu }^{1}(s) &=&\frac{\left( w-1\right) \sin [\left( w+1\right) t]}{%
2w\left( w+1\right) }+\frac{\left( w+1\right) \sin [\left( w-1\right) t]}{%
2w\left( w-1\right) } \\
\gamma _{\mu }^{2}(s) &=&\frac{\left( w+1\right) \cos [\left( w-1\right) t]}{%
2w\left( w-1\right) }+\frac{\left( w-1\right) \cos [\left( w+1\right) t]}{%
2w\left( w+1\right) } \\
\gamma _{\mu }^{3}(s) &=&-\frac{\cos \left( s\right) }{\mu w}.
\end{eqnarray*}%
and $w=\frac{\sqrt{1+\mu ^{2}}}{\mu },$ \ $\mu \in IR-\left\{ 0\right\} .$
By using the theorem \ref{c5}, we can obtained the general equation of
general helix family in Euclidean 3-space according to the non-zero constant 
$\mu $ as follows%
\begin{equation*}
\widetilde{\gamma }_{\mu }\left( s_{\ast }\right) =\left( \widetilde{\gamma }%
_{\mu }^{1}(s),\widetilde{\gamma }_{\mu }^{2}(s),\widetilde{\gamma }_{\mu
}^{3}(s)\right)
\end{equation*}%
where%
\begin{eqnarray*}
\widetilde{\gamma }_{\mu }^{1}(s) &=&\frac{\left( c-s\right) }{2w}\left\{
\left( w+1\right) \cos [\left( w-1\right) t]+\left( w-1\right) \cos [\left(
w+1\right) t]\right\} \\
&&+\frac{1}{2w\left( w^{2}-1\right) }\left\{ \left( w+1\right) ^{2}\sin
[\left( w-1\right) t]+\left( w-1\right) ^{2}\sin [\left( w+1\right)
t]\right\} \\
\widetilde{\gamma }_{\mu }^{2}(s) &=&\frac{-\left( c-s\right) }{2w}\left\{
\left( w+1\right) \sin [\left( w-1\right) t]+\left( w-1\right) \sin [\left(
w+1\right) t]\right\} \\
&&+\frac{1}{2w\left( w^{2}-1\right) }\left\{ \left( w+1\right) ^{2}\cos
[\left( w-1\right) t]+\left( w-1\right) ^{2}\cos [\left( w+1\right)
t]\right\} \\
\widetilde{\gamma }_{\mu }^{3}(s) &=&\frac{\left( c-s\right) \sin \left(
s\right) -\cos \left( s\right) }{\mu w}
\end{eqnarray*}%
and $w=\frac{\sqrt{1+\mu ^{2}}}{\mu },$ \ $\mu \in IR-\left\{ 0\right\} ,$
with the parametrization%
\begin{equation*}
\frac{ds}{ds_{\ast }}=\frac{\mu w\sqrt{2}}{\left( c-s\right) \sqrt{\left(
\mu ^{2}\left( w^{2}-1\right) ^{2}+1\right) \left( \cos \left( 2t\right)
+1\right) }}
\end{equation*}%
and the curvature and torsion of curve $\widetilde{\gamma }_{\mu }\left(
s_{\ast }\right) $ is 
\begin{equation*}
\widetilde{\kappa }=\frac{\sqrt{2}}{\left( c-s\right) \sqrt{\left( \cos
\left( 2t\right) +1\right) }}\text{ \ , \ }\widetilde{\tau }=\frac{-\mu }{%
\left( c-s\right) \cos \left( t\right) }.
\end{equation*}
\begin{figure}[htbp]
\label{slanthelix} \centering
\includegraphics[width=3cm,angle=0,height=5cm]{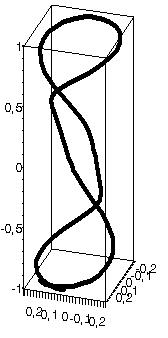}
\caption{Slant helix}
\end{figure}
\begin{figure}[htbp]
\label{generalhelix} \centering
\includegraphics[width=3cm,angle=0,height=5cm]{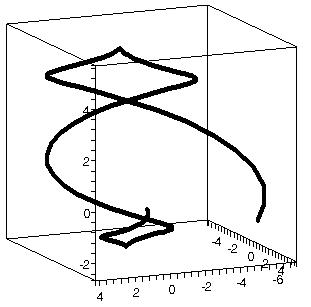}
\caption{General helix with $c=-1, {\protect\mu}=\pm 1/4$.}
\end{figure}
\begin{figure}[htbp]
\label{tangentindicatrix} \centering
\includegraphics[width=3cm,angle=0,height=3cm]{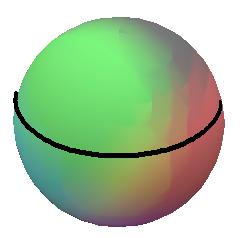}
\caption{Tangent indicatrix of involute of a slant helix.}
\end{figure}
\begin{figure}[htbp]
\label{principalindicatrix} \centering
\includegraphics[width=3cm,angle=0,height=3cm]{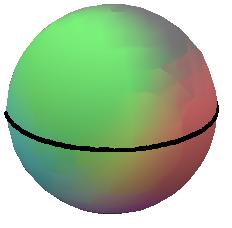}
\caption{Principal normal indicatrix of involute of a slant helix.}
\end{figure}
\begin{figure}[htbp]
\label{binormalindicatrix} \centering
\includegraphics[width=3cm,angle=0,height=3cm]{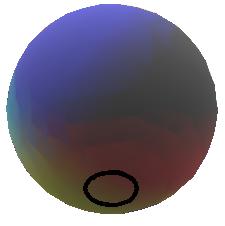}
\caption{Binormal indicatrix of involute of a slant helix.}
\end{figure}
\end{example}

\end{document}